\documentclass[review,authoryear]{elsarticle}
\usepackage{lineno}
\modulolinenumbers[5]
\usepackage{array,epsfig,fancyheadings,rotating}
\usepackage{sectsty,secdot}
\usepackage{amsmath}
\usepackage{amssymb}
\usepackage{amsfonts}
\usepackage{colortbl}
\usepackage{arydshln}
\usepackage{multirow}
\usepackage{multicol}
\usepackage{amsthm}
\usepackage{amssymb}
\usepackage{tipa}
\usepackage{graphicx,color,overpic}
\usepackage{lscape}
\usepackage{algorithmic,algorithm}
\usepackage{epstopdf,epsfig}
\usepackage{booktabs}
\usepackage{natbib}
\usepackage{ulem}
\newtheorem{theorem}{Theorem}
\newtheorem{lemma}{Lemma}
\newtheorem{corollary}{Corollary}

\theoremstyle{definition}

\newtheorem{example}{Example}
\newtheorem{remark}{Remark}

\newcommand{\ba}{\begin{array}}
\newcommand{\ea}{\end{array}}
\newcommand{\bt}{\begin{tabular}}
\newcommand{\et}{\end{tabular}}
\newcommand{\btb}{\begin{table}}
\newcommand{\etb}{\end{table}}
\newcommand{\bc}{\begin{center}}
\newcommand{\ec}{\end{center}}
\newcommand{\bea}{\begin{eqnarray}}
\newcommand{\eea}{\end{eqnarray}}
\newcommand{\Bea}{\begin{eqnarray*}}
\newcommand{\Eea}{\end{eqnarray*}}
\newcommand{\beq}{\begin{equation}}
\newcommand{\eeq}{\end{equation}}
\newcommand{\bes}{\begin{subequations}}
\newcommand{\ees}{\end{subequations}}
\newcommand{\bal}{\begin{aligned}}
\newcommand{\eal}{\end{aligned}}

\def \bfm#1{\mbox{\boldmath$#1$}}

\def \A {{\bfm A}}   
  
 \def \D {{\bfm D}}

  \def \0 {{\bfm 0}}

\def \s {{\bfm s}} 
\def \t {{\bfm t}}

\def \x {{\bfm x}} \def \X {{\bfm X}}
 \def \y {{\bfm y}}
 \def \z {{\bfm z}}

\def \one {{\bf 1}} \def \zero {{\bf 0}}

 \def \cd {{\mathcal D}} 
 \def \ck {{\mathcal K}}

\def \cu {{\mathcal U}}

\makeatletter

\newcommand{\Rmnum}[1]{\expandafter\@slowromancap\romannumeral #1@}
\makeatother

\newcommand{\revA}[1]{{#1}}

\def \d {{\rm d}}

\journal{Statistics }
\begin{document}

\begin{frontmatter}

\title{Uniformity criterion for designs with both qualitative and quantitative factors}

\author{Mei Zhang$^{1}$}


\author{Feng Yang$^{2}$}

\author{Yong-Dao Zhou$^3$\corref{mycorrespondingauthor}}
\cortext[mycorrespondingauthor]{Corresponding author}
\ead{ydzhou@nankai.edu.cn}

\address{1.College of Mathematics, Sichuan University, Chengdu 610064, China}
\address{2.School of Mathematical Science, Sichuan Normal University, Chengdu 610068, China }
\address{3.School of Statistics and Data Science, LPMC $\&$ KLMDASR, Nankai University, Tianjin 300071, China }

\begin{abstract}
Experiments with both qualitative and quantitative factors occur frequently in practical applications. Many construction methods for this kind of designs, such as marginally coupled designs, were proposed to pursue some good space-filling structures. However, few criteria can be adapted to quantify the space-filling property of designs involving both qualitative and quantitative  factors. As the uniformity is an important space-filling property of a design, in this paper, a new uniformity criterion, qualitative-quantitative discrepancy (QQD), is proposed for assessing the uniformity of designs with both types of factors. The closed form and lower bounds of the QQD are presented to calculate the exact QQD values of designs and recognize the uniform designs directly. In addition, a connection between the QQD and the balance pattern is derived, which not only helps to obtain a new lower bound but also provides a statistical justification of the QQD. Several examples show that the proposed criterion is reasonable and useful since it can distinguish distinct designs very well.
\end{abstract}


\end{frontmatter}


\section{Introduction} \label{intro}

Experimental designs with qualitative and quantitative factors
have received growing attention in recent years,
see \cite{H09}, \cite{QWW08} and \cite{Z11}
for \revA{computer experiments}; \cite{WD98}, \cite{A00}, \cite{C03} and \cite{TB03} for  response surface designs.
For computer experiments, there are two systematic approaches to handle this issue.
\cite{Q12} suggested using sliced Latin hypercube designs
to accommodate qualitative and quantitative factors,
where each slice of quantitative factors corresponds to
a level combination of qualitative factors.
However, the run sizes of such designs can be very large even
for a moderate number of qualitative factors.
Inspired by this, \cite{DHL15} proposed a new class of designs,
marginally coupled designs (MCDs), of which the design points
for quantitative factors form a Latin hypercube design (LHD),
and for each level of any qualitative factor,
the corresponding design points for quantitative factors
compose a smaller LHD. Intuitively, MCDs have
some desirable space-filling properties.
A series of construction methods for MCDs
had been presented,
see \cite{DHL15}, \cite{HLS17}  and \cite{HLS19}.
Except for \revA{designs for computer experiments with both qualitative and quantitative factors}, the
response surface designs with qualitative and
quantitative factors have been studied and applied
extensively in science and engineering applications.
\cite{WD98} proposed a general approach for constructing
response surface designs of economical size with
qualitative and quantitative factors. \cite{A00} applied a
dual response surface optimization technique for constructing
robust response surface designs with quantitative and
qualitative factors. \cite{TB03} suggested a methodology for
developing a simulation metamodel involving both quantitative and qualitative factors to deal with various strategic issues, such as metamodel estimation, analysis, comparison, and validation. Moreover, the response surface methodology was also developed in \cite{C14}, to investigate the photocatalytic degradation of phenol and phenol derivatives using a Nano-TiO$_2$ catalyst.
In their study, two quantitative factors (TiO$_2$ particle size and temperature) and one qualitative factor (reactant type) were considered.

\begin{figure}[t!]
\centering
\includegraphics[width=0.7\textwidth]{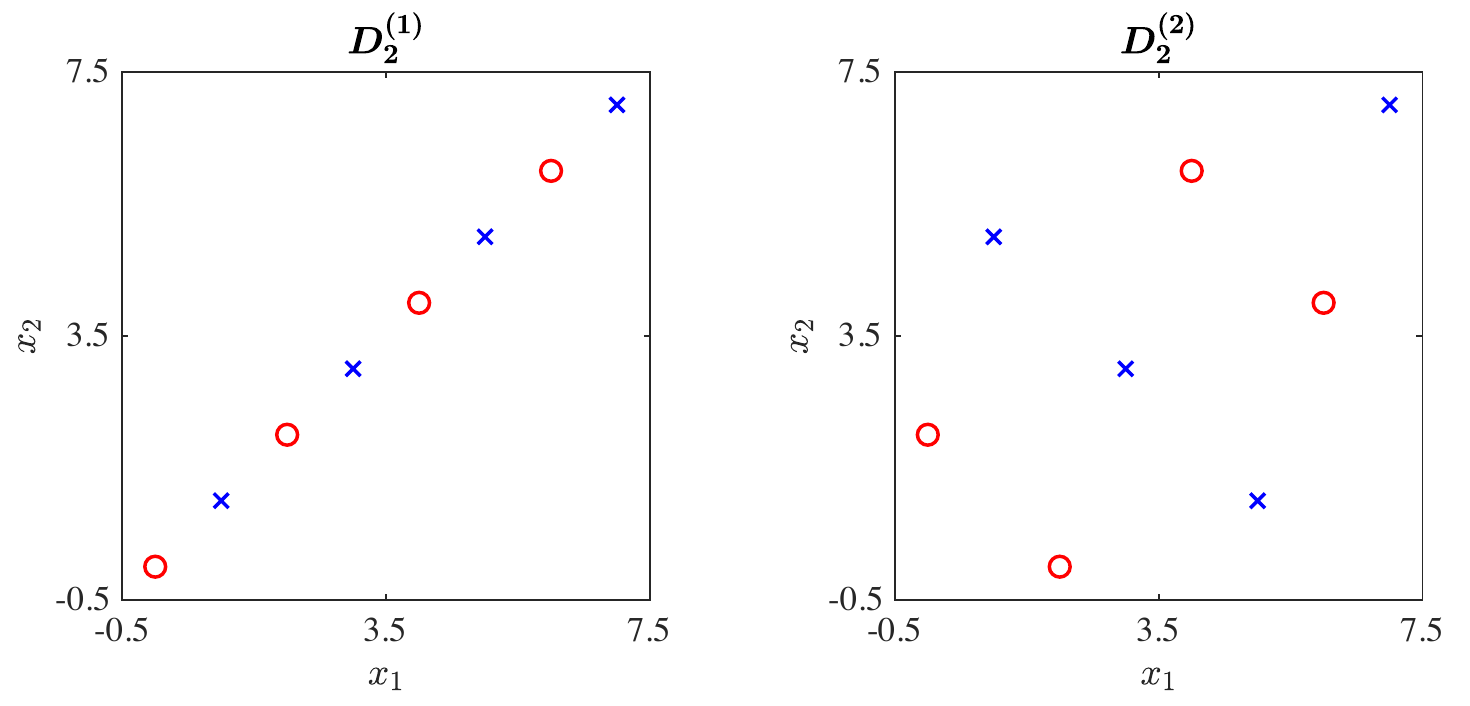}
\caption{The scatter plots for  the
$\D_2^{(1)}$ and $\D_2^{(2)}$  in
(\ref{eg:1+2}), where the points of $\D_2^{(j)} ( j = 1, 2 )$
corresponding to level 0 and 1 of $\D_1$ are represented
by ``\revA{$\circ$}" and ``\revA{\footnotesize{$\times$}}", respectively.}
\label{Fig:1+2}
\end{figure}

Let $\D=(\D_1,\D_2)$ be a design with $n$ runs for $p$
qualitative and $q$ quantitative factors,
in which $\D_1$ and $\D_2$ are sub-designs for
qualitative and quantitative factors, respectively.
Given the parameters $(n,p,q)$, a great number of
MCDs can be constructed.
It arises a natural issue that: which one is the best.
In other words, it needs some proper criterion to distinguish
these MCDs. For example, consider two MCDs
with $8$ runs for one qualitative and two quantitative factors,
\bea\label{eg:1+2}\D^{(1)} = (\D_1,\D_2^{(1)})
\text{~~~and~~~} \D^{(2)} = (\D_1,\D_2^{(2)})\eea
where $\D_1 = ( 0, 0, 0, 0, 1, 1, 1, 1 )^T$,
$\D_2^{(1)} =\big((0, 2, 4, 6, 1, 3, 5, 7)^T,
(0, 2, 4, 6, 1, 3, 5, 7)^T\big)$, and
$\D_2^{(2)} =\big((0, 2, 4, 6, 1, 3, 5, 7)^T,
(2, 0, 6, 4, 5, 3, 7, 1)^T\big)$.
Figure \ref{Fig:1+2} shows the scatter plots of sub-designs
$\D_2^{(1)}$ and $\D_2^{(2)}$ with respect to the levels of $\D_1$.
From the intuition, $\D^{(2)}$ outperforms $\D^{(1)}$
since the whole points of $\D_2^{(1)}$
as well as the corresponding points of $\D_2^{(1)}$
with respect to each level of $\D_1$,
locate on the diagonal line, respectively.

If the number of factors is relatively small, visually,
the scatter plots can help us to make a judgment,
whereas for the high-dimensional cases (even for $p>1$ and $q>2$),
it may be impractical.
Motivated by this,
a space-filling criterion is demanded to compare different designs
with qualitative and quantitative factors
especially for high-dimensional cases.

\revA{In computer experiments, the maximin distance criterion proposed by \cite{J90} is a popular used criterion for measuring the space-filling property of designs.
However, it only considers the space-filling property on the full-dimensional space and may result in poor projections onto lower-dimensional spaces.  \cite{J15} proposed the maximum projection (MaxPro) criterion, which can simultaneously optimize the space-filling properties of the design points with respect to all possible subsets of factors. And the MaxPro criterion is  suitable for the cases where any level is distinct from one another in each factor.
To handle more cases in computer experiments, \cite{J19} extends the MaxPro criterion to accommodate multiple types of factors. 
From another version,} the uniformity is another important criterion for assessing the
space-filling property of a design,
see \cite{FW94} and \cite{FLS06}.
As the measure of the uniformity of designs,
discrepancies have played a key role in the theory of uniform designs.
So far, many discrepancies have been proposed to
measure the uniformity of designs on continuous
design space $[0,1]^m$.
For example, \cite{H98,H99} used the tool of reproducing kernel Hilbert
spaces to define the centered $L_2$-discrepancy (CD) and
the wrap-around $L_2$-discrepancy (WD), and
\cite{ZFN13} proposed the mixture discrepancy (MD).
Nevertheless, these discrepancies are only applicable for designs
with quantitative factors.

For designs with a finite number of levels, similar to
CD and WD, \cite{HL02} and \cite{LH02}
proposed the discrete discrepancy (DD) that was directly
defined on a discrete domain. It seems that
DD is the only discrepancy we can use when the inputs
are qualitative. Certainly, if we recognize \revA{that} the levels of quantitative
inputs are taken from a finite set, the DD also works for
quantitative factors.
Hence, the DD could be geared to the designs
with both types of factors.
While \cite{ZNS08} pointed out that the DD
is not an effective uniformity criterion for designs
with multi-level quantitative factors as
the DD is based on the Hamming distance.
For instance, it can be easily checked that
$\D^{(1)}$ and $\D^{(2)}$ in (\ref{eg:1+2})
have the same DD-value, which also reaches the lower
bound of DD in \cite{QF04}. That is,
$\D^{(1)}$ and $\D^{(2)}$ are the best designs under the DD,
simultaneously. This result is not in accord with the intuition
since $\D^{(2)}$ should be more uniform than $\D^{(1)}$
from Figure \ref{Fig:1+2}, which indicates that the DD is not suitable for
assessing the uniformity of the designs
with both qualitative and multi-level quantitative factors.
Therefore, it calls for a new uniformity criterion for
designs containing both quantitative and qualitative factors.
To address this issue, in this paper, we propose
a new uniformity criterion
\textit{qualitative-quantitative discrepancy},
to assess the uniformity of designs with both types of factors.

The remainder of this paper is organized as follows.
Section \ref{se_qd} first uses the Hickernell's approach
to derive the proposed discrepancy, QQD, then
gives the quadratic form of the QQD
for finding the uniform designs.
The lower bounds of the QQD
that provides a simple way to recognize the uniform designs,
are derived in Section \ref{se_lb}. Section \ref{relationships}
shows the relationships between QQD and balance pattern,
and induces another lower bound.
Some illustrative examples are given in Section \ref{egs}
to explain the reasonability and application of QQD,
and demonstrate the tightness of the obtained
lower bounds in Sections \ref{se_lb} and \ref{relationships}.
The last section gives some conclusions. All the proofs are listed in the Appendix.

\section{The qualitative-quantitative discrepancy}\label{se_qd}

In the field of uniform 
\revA{designs, U-type designs are generally preferred}
because of their projection uniformity on each dimension.
A  U-type design  is an $n\times m$ matrix
$\X=(\d_1,\dots,\d_m)$, of which the $k$th column
takes values from the set of $\{0,1,\dots,s_k-1\}$ equally often.
Denote all of such designs by $\cu(n,s_1\cdots s_{m})$.
If some $s_k$'s are equal, we rewrite them as
$\cu(n,s_1^{m_1} \cdots  s_l^{m_l})$ with
$m=\sum_{k=1}^l m_k$.
Denote $\cu(n,s_1\cdots s_p s_{p+1} \cdots s_{p+q})$
as all
the U-type designs whose first $p$ factors are qualitative
and the last $q$ factors are quantitative, then
$\D=(\D_1,\D_2) \in \cu(n,s_1\cdots s_p s_{p+1} \cdots s_{p+q})$
implies that
$\D_1\in\cu(n,s_1\cdots s_p)$ and $\D_2\in
\cu(n,s_{p+1}\cdots s_{p+q})$.
Specially, when $s_{p+1} = \cdots  = s_{p+q} = n$, $\D_2$ is an LHD($n,q$). \revA{According to the structure of MCDs, MCDs are also U-type designs.}

We first review the derivations of
the discrepancies, and give the formula of the QQD.
Then, the quadratic form for the QQD
is induced to help to obtain the uniform design under the
QQD.

\subsection{The expression of the QQD}
Let $F$ be the uniform distribution function on the design
region $\chi~(=\chi_1\times\dots\times\chi_{m})$,
and $F_n$ be the
empirical distribution of a set of points
$\X=\{\x_1,\dots,\x_n\}$.
For a given kernel function $\ck(\t,\z)$, any discrepancy of
design $\X$ with $n$ runs and $m$ factors can be defined
by (see \cite{F18})
\begin{align}\label{discre}
\cd^2(\X,\ck)&= \int_{\chi^{2}}\ck(\t,\z)\mathrm{d}(F(\t)-F_n(\t))
\mathrm{d}(F(\z)-F_n(\z)) \nonumber\\
&=\int_{\chi^{2}}\ck(\t,\z)\mathrm{d}F(\t)
\mathrm{d}F(\z)-\frac{2}{n}\sum_{i=1}^{n}
\int_{\chi}\ck(\t,\x_i)\mathrm{d}F(\t)+\frac{1}{n^2}\sum_{i,j=1}^{n}
\ck(\x_i,\x_j).
\end{align}
The discrepancy in (\ref{discre}) is devoted to measuring
the difference between $F_n$ and $F$, hence,
the smaller its value is,
the more uniformly the design points spread on $\chi$.
A design is called a \textit{uniform design} if it has the
minimal discrepancy over the design space.
From (\ref{discre}),
the discrepancy is determined by the
kernel function $\ck(\cdot,\cdot)$, that is, a kernel
function will give rise to a type of discrepancy.
Usually, $\ck(\cdot,\cdot)$ has a product form
\begin{align}\label{kernel}
\mathcal{K}(\t,\z)= \prod_{k=1}^{m}\mathcal{K}_k(t_k,z_k),
\end{align}
where $\mathcal{K}_k(t_k,z_k)$ represents the kernel
of the $k$th factor.
It should be noted that the $\cd^2(\X,\ck)$ defined in (\ref{discre}) with the product formed kernel $\ck(\cdot,\cdot)$ in (\ref{kernel}) considers the space-filling properties on projections to all subsets of factors.
For any subset $\revA{u} \subseteq \{1:m\}$, let $\chi_{\revA{u}}$ be the experimental region of the corresponding factors whose indices are in $\revA{u}$, and $\t_{\revA{u}}$ be the vector containing the components of $\t$ indexed by $\revA{u}$.
Denote $\widetilde{\mathcal{K}}_{\revA{u}}(\t_{\revA{u}},\z_{\revA{u}}) = \prod_{k\in \revA{u}}\left[\mathcal{K}_k(t_k,z_k) - 1\right]$ when $\revA{u} \neq \emptyset$ and $1$ otherwise. Then $\mathcal{K}(\t,\z) = \sum_{\revA{u} \subseteq \{1:m\}} \widetilde{\mathcal{K}}_{\revA{u}}(\t_{\revA{u}},\z_{\revA{u}})$ by the induction. At this moment, the discrepancy in (\ref{discre}) could be rewritten as
\begin{align*} 
\cd^2(\X,\ck)&= \int_{\chi^{2}}\sum_{\revA{u} \subseteq \{1:m\}} \widetilde{\mathcal{K}}_{\revA{u}}(\t_{\revA{u}},\z_{\revA{u}})\mathrm{d}(F(\t)-F_n(\t))\mathrm{d}(F(\z)-F_n(\z)) \nonumber\\
&=\sum_{\revA{u} \subseteq \{1:m\}}\int_{\chi_{\revA{u}}^{2}}\widetilde{\mathcal{K}}_{\revA{u}}(\t_{\revA{u}},\z_{\revA{u}})\mathrm{d}(F(\t)-F_n(\t))\mathrm{d}(F(\z)-F_n(\z)).
\end{align*}
\revA{It} implies that the defined $\cd^2(\X,\ck)$ could measure not only the uniformity of $\X$ on $\chi$, but also projection uniformity of $\X$ on any $\chi_{\revA{u}}$, where ${\revA{u}}$ is a non-empty subset of $\{1, \dots, m\}$.

Recall that
the WD is valid for the design with continuous inputs.
If one transforms its levels into the unit cube $[0,1]$ by
$(2x+1)/(2s_k),x=0,\dots,s_k-1$,
the kernel function for the $k$th factor is defined by
$\mathcal{K}_k(t_k,z_k)=\frac{3}{2}-|~t_k-z_k~|+|~t_k-z_k~|^2$.
In addition, the DD
is intended for the designs with finite numbers of levels,
whose kernel function of the $k$th factor is
$\mathcal{K}_k(t_k,z_k)=a^{\delta_{t_kz_k}}b^{1-\delta_{t_kz_k}}$
with $a>b$, where $\delta_{t_kz_k}$ equals $1$ if $t_k=z_k$ and
$0$ otherwise, and the design region is
$\chi_k=\{0,1,\dots,s_k-1\}$.
For a design $\D=(\D_1,\D_2)\in\cu(n,s_1\cdots s_p
s_{p+1} \cdots s_{p+q})$.
It is reasonable to allocate the kernel of DD to the
qualitative factors of $\D$ and the kernel of WD to the
quantitative factors of $\D$. By doing this,
the kernel function for the $k$th factor of $\D$ is
\begin{align}\label{ke_qd}
\mathcal{K}_k(t_k,z_k)=
\begin{cases}
~a^{\delta_{t_kz_k}}b^{1-\delta_{t_kz_k}}, & \text{for $k=1,\dots,p$,}\\
~\frac{3}{2}-|~t_k-z_k~|+|~t_k-z_k~|^2, & \text{for $k=p+1,\dots,p+q$.}
\end{cases}
\end{align}
We call the discrepancy with respect to kernel function in
(\ref{ke_qd}) as the qualitative-quantitative discrepancy.
Based on (\ref{discre}), (\ref{kernel}) and (\ref{ke_qd}),
it can be seen that the value of QQD
depends on the parameters $a$ and $b$.
Thus the choice of $a$ and $b$ is a crucial issue
for the proposed QQD. Let us focus on the kernel functions for the two types of factors.
The kernel function
for the qualitative factors, $ \ck_k(x_{ik},x_{jk}),k=1,\dots,p,$
has  the lower bound $b$ and upper bound $a$,
while the kernel function
for the quantitative factors,
$
\ck_k(x_{ik},x_{jk}),k=p+1,\dots,p+q
$,
has the lower bound $5/4$ and the upper bound $3/2$.
In general, if we do not have any prior information about
the importance of factors, just treat them comparably.
In this sense, let $a=3/2$ and $b=5/4$ in (\ref{ke_qd}),
so that the kernel functions for both types of
factors have the same ranges, which indicates that the
qualitative and quantitative factors dominate
the QQD equally. In the rest of the paper, unless
otherwise specified, we take $a=3/2$ and $b=5/4$.
Under this, we can derive the closed form of the expression for QQD.

\begin{theorem}\label{th_qd} For a design $\D\in\cu(n,s_1\cdots s_ps_{p+1} \cdots s_{p+q})$,
the expression of the squared QQD is as follows
\begin{align}
\text{QQD}^2(\D)=C+\frac{1}{n^2}\sum_{i,j=1}^{n}\left(\frac{5}{4}\right)^p
\left(\frac{6}{5}\right)^{\delta_{ij}(\D_1)}\cdot\prod_{k=p+1}^{p+q}
\left(\frac{3}{2}-|~x_{ik}-x_{jk}~|+|~x_{ik}-x_{jk}~|^2\right).\label{qd_ex}
\end{align}
where
$C = -\prod_{k=1}^{p}\left(\frac{5s_k+1}{4s_k}\right)\left(\frac{4}{3}\right)^q,$
$\delta_{ij}(\D_1)$ represents the coincidence number
between the $i$th and $j$th rows of $\D_1$.
\end{theorem}

It is noted that, although the proposed QQD is motivated by comparing different MCDs, it is not limited to MCDs. According to the expression in (\ref{qd_ex}), QQD can be used for quantifying the uniformity of any design with both qualitative and quantitative factors.

We now review the MCDs in (\ref{eg:1+2})
to illustrate that the
proposed criterion, QQD, is more effective and appropriate
than DD to measure the goodness of different designs
with both types of factors.
\begin{example} Based on the expression of QQD in
Theorem \ref{th_qd},  for the two MCDs in (\ref{eg:1+2}),
$\text{QQD}^2(\D^{(1)}) = 0.0213$ and
$\text{QQD}^2(\D^{(2)}) = 0.0164$, respectively.
That is to say under the QQD criterion,
design $\D^{(2)}$ is better than design $\D^{(1)}$, which
 is coincident with the \revA{scatter plot} in Figure \ref{Fig:1+2}.
 As discussed in   
 \revA{Section \ref{intro}}, the DD cannot distinguish
 the two different designs at all.
Consequently, the QQD is more sensitive than the DD to
evaluate the uniformity of designs associated
with qualitative and quantitative factors.
\revA{Here, we also compute the extended MaxPro criterion value of $\D^{(1)}$ and $\D^{(2)}$ by the expression in \cite{J19} and obtain the corresponding values $8.8320$ and $6.3366$. It means design $\D^{(2)}$ is better than design $\D^{(1)}$ under this criterion. The agreement with the result under the QQD criterion further certifies the reasonability of QQD.}
\end{example}

\revA{In Example 1, it is reflected that both the QQD and the  extended MaxPro criterion perform well when they are used to measure the space-filling property of designs with both qualitative and quantitative factors, 
although the two criteria are defined from different points.} 

\subsection{The quadratic form of the QQD}
In this subsection, the uniform designs among the design
space under the QQD criterion are discussed.
We first present the quadratic form of the QQD
and then solve the optimization problem of a special
convex quadratic programming to obtain the uniform designs.

For ease of presentation, given a design
$\D=(\D_1,\D_2)\in \cu(n,s_1\cdots s_p s_{p+1} \cdots s_{p+q})$,
a new notation is introduced. Let $\y=\y(\D)$ be an
$N~( = \prod_{k = 1}^{p+q}s_k)$-dimensional column vector
with components $n(i_1,\dots,i_{p+q})$ arranged
lexicographically, where $n(i_1,\dots,i_{p+q})$ is the
number of runs at the level combination $(i_1,\dots,i_{p+q})
$ in design $\D$. Obviously, $\y$ is a non-negative vector and satisfies $\one_N^T\y = n$.
Moreover, if each element of $\y$ is an integer,
the design is called an exact design, otherwise,
a continuous design. Any U-type design is an exact design.
Based on the concept of $\y(\D)$,
the expression of $\text{QQD}^2(\D)$ in (\ref{qd_ex})
can be rewritten as a quadratic form given in the following lemma.
\begin{lemma} \label{quadratic_form}
If $\D \in \cu(n,s_1\cdots s_{p+q})$
and $\y=\y(\D)$,
we have \begin{align*}
\text{QQD}^2(\D) = -\prod_{k=1}^{p}\left(\frac{5s_k+1}{4s_k}\right)
\left(\frac{4}{3}\right)^q+\frac{1}{n^2}\y^T\A\y,
\end{align*}
where $\A=\A_1\bigotimes \A_2 \bigotimes \cdots \bigotimes \A_{p+q}$,
$\A_k=(t_{ij}^k), i, j=1,\dots, s_k, k=1,\dots, p+q$,
and $\bigotimes$ is the Kronecker product, furthermore, \bc
$t_{ij}^k= \begin{cases}
~\left(\frac{3}{2}\right)^{\delta_{ij}}\left(\frac{5}{4}\right)^{1-\delta_{ij}}, & \text{for $k=1,\dots,p$},\\
~\frac{3}{2}-\frac{|i-j|(s_k-|i-j|)}{s_k^2}, &
\text{for $k=p+1,\dots,p+q$.}
\end{cases}$\ec
\end{lemma}

Lemma \ref{quadratic_form} indicates that,
essentially, $\text{QQD}^2(\D)$ is a convex function of $\y$ with
some constraints. Provided this quadratic form, the
theory of convex optimization can be utilized
to solve the optimization problem of
minimizing $\text{QQD}^2(\D)$.
The corresponding solution $\D^*$ is the
uniform design under the QQD criterion.
The rest of this section aims at exploring
the uniform designs under the QQD.
The proof of Theorem \ref{th_uni} is analogous to that
of Theorem 1 in \cite{ZFN12}, which relies on Lemma
\ref{le_property} in Appendix. For simplicity, we omit it.

\begin{theorem}\label{th_uni}
The design $\D^*$ minimizes $\text{QQD}^2(\D)$ over
$\cu(n,s_1 \cdots s_{p+q})$ when
$\y(\D^*)=\frac{n}{N}\one_N$.
\end{theorem}

Theorem \ref{th_uni} provides a sufficient condition for
a design $\D\in \cu(n,s_1\cdots s_{p+q})$
being a uniform design under QQD.
It should be noted that,
in the proof of Theorem \ref{th_uni}, the constrains of $\y$ are relaxed to $\y \geq \zero_N$ and $\one_{N}^T\y = n$ as in \cite{ZFN12}. Thus the obtained optimal design
$\D^*$ under the QQD is not an exact design unless
$n ~(\text{mod}~N) =0$.
\revA{T}he uniform design $\D^*$ with $\y(\D^*)=\frac{n}{N}\one_N$ implies that if the points at each level combination have the same weight in a design, the design achieves the best uniformity under \revA{the} QQD criterion.
For an exact uniform design, the squared QQD can be
calculated by a concise expression, of parameters only,
given by the following Corollary \ref{co_full}.

\begin{corollary} \label{co_full}
When $n=cN$ with $c$ being a positive integer,
the $\D^*$ in Theorem~\ref{th_uni}
is a repetition of a full factorial design and its
squared QQD is
\begin{align*}
\text{QQD}^2(\D^*) =-\prod_{k=1}^{p}\left(\frac{5s_k+1}{4s_k}\right)
\left(\frac{4}{3}\right)^q + \prod_{k=1}^{p}\left(\frac{5s_k+1}{4s_k}\right)
\prod_{k=p+1}^{p+q}\left(\frac{4}{3}+
\frac{1}{6s_k^2}\right).
\end{align*}
\end{corollary}

 Corollary \ref{co_full} shows that
all of the  full factorial designs and their repetitions
are uniform designs under the QQD,
since the formula of QQD in Corollary \ref{co_full}
is independent of the repetition times $c$.
This result \revA{is natural}
since full factorial designs are the best designs under many existing criteria, such as WD and minimum aberration.

\section{The lower bounds of QQD}\label{se_lb}

The lower bounds of a discrepancy can be employed as
a benchmark not only in searching for uniform designs but also
in helping validate that some good designs are uniform.
The issue of lower bounds for all kinds of discrepancies
have been considered and
many researchers have made plenty of
effort in finding the lower bounds for different discrepancies,
such as CD, WD, MD and DD,
see \cite{F18} for a comprehensive review.
In this section, some lower bounds of the QQD for U-type designs
with qualitative and quantitative factors are explored.

\begin{theorem}\label{th_lbqd} Let $\D=(\D_1,\D_2)$ with
$\D_1\in \cu(n,s_1\cdots s_p)$ and
$\D_2\in \cu(n,s_{p+1} \cdots s_{p+q})$,
and assume $s_{p+1},\dots,s_{t}$ be odd and
$s_{t+1},\dots,s_{p+q}$ be even, where $p\leq t \leq p+q$,
then $QQD^2(\D)\geq \text{LB}_1$, where
\begin{align*} \text{LB}_1
= &
~C
+\frac{1}{n}\left(\frac{3}{2}\right)^{p+q}
+\frac{(n-1)}{n}\left(\frac{5}{4}\right)^p\left(\frac{6}{5}\right)^
{\sum_{k=1}^{p}\frac{n-s_k}{s_k(n-1)}}
\left(\frac{3}{2}\right)^{\sum_{k=p+1}^{p+q}\frac{n-s_k}{s_k(n-1)}}
\left(\frac{5}{4}\right)^{\sum_{k=t+1}^{p+q}\frac{n}{s_k(n-1)}}
\nonumber\\
&\times \prod_{k=p+1}^{t}\prod_{i=1}^{(s_k-1)/2}\left(\frac{3}{2}-
\frac{2i(2s_k-2i)}{4s_k^2}\right)^{\frac{2n}{s_k(n-1)}}
\prod_{k=t+1}^{p+q}\prod_{i=1}^{(s_k/2)-1}\left(\frac{3}{2}-
\frac{2i(2s_k-2i)}{4s_k^2}\right)^{\frac{2n}{s_k(n-1)}},
\end{align*}
where $C$ is given in Theorem \ref{th_qd}.
\end{theorem}

Theorem \ref{th_lbqd} obtains the results for general
cases, i.e.,
any U-type design $\D$ with both qualitative and
quantitative factors\revA{,
where the number of quantitative factors with odd levels could be $0, 1, \dots, q.$}
Specially, $t = p + q$ means the levels of quantitative factors are all odd and $t = p$ means the levels of quantitative factors are all even.
This provides a benchmark for identifying the uniform designs.
For instance,
if the QQD of a design \revA{reaches}
the corresponding lower bound, such \revA{a} design must be a uniform design.
Under some special parameter settings,
the associated results can be directly acquired and described in
Corollary \ref{co_sy} and Remark \ref{re_redu}.

\begin{corollary}\label{co_sy}
Let $\D=(\D_1,\D_2)\in \cu(n,s_1^ps_2^q)$,
where $\D_1$ and $\D_2$ are
symmetric designs with $s_1$ levels and $s_2$
levels, respectively,
then the lower bound of $QQD^2(\D)$ is given by
\begin{align*}
LB_{odd}
 =  - & \left(\frac{5s_1+1}{4s_1}\right)^p
\left(\frac{4}{3}\right)^q+\frac{1}{n}\left(\frac{3}{2}\right)^{p+q}
+\frac{(n-1)}{n}\left(\frac{5}{4}\right)^p\left(\frac{6}{5}\right)^{\frac{p(n-s_1)}{s_1(n-1)}}
\left(\frac{3}{2}\right)^{\frac{q(n-s_2)}{s_2(n-1)}}\nonumber\\
&\times\prod_{i=1}^{(s_2-1)/2}\left(\frac{3}{2}
-\frac{2i(2s_2-2i)}{4(s_2)^2}\right)^{\frac{2nq}{s_2(n-1)}},
\text{~~~~~~~~~~~~~~~~~~~~~~for odd $s_2$,}\end{align*}
and \begin{align*}
LB_{even}
=  - & \left(\frac{5s_1+1}{4s_1}\right)^p
\left(\frac{4}{3}\right)^q+\frac{1}{n}\left(\frac{3}{2}\right)^{p+q}
+\frac{(n-1)}{n}\left(\frac{5}{4}\right)^p\left(\frac{6}{5}\right)^{\frac{p(n-s_1)}{s_1(n-1)}}
\left(\frac{3}{2}\right)^{\frac{q(n-s_2)}{s_2(n-1)}}\nonumber\\
&\times\left(\frac{5}{4}\right)^{\frac{nq}{s_2(n-1)}}
\prod_{i=1}^{(s_2/2)-1}\left(\frac{3}{2}-
\frac{2i(2s_2-2i)}{4(s_2)^2}\right)^{\frac{2nq}{s_2(n-1)}},
\text{~~~~~~for even $s_2$.}
\end{align*}
\end{corollary}

\begin{remark}\label{re_redu} If $q=0$, the QQD reduces to DD
in \cite{HL02, LH02} and the lower bound in Theorem \ref{th_lbqd}
 becomes that for DD of Theorem 2 in \cite{FLL03}.
When $p=0$, the QQD becomes WD in \cite{H98} and
the lower bound for QQD becomes that for WD of Theorem
1 in \cite{ZN08}.
\end{remark}


\section{The connection between QQD and balance pattern}\label{relationships}

The concept of balance pattern was introduced by
\cite{FLW03},
which characterizes the column balance of a design.
For a two-level or three-level U-type design, \cite{FLW03}
pointed out that  WD can be expressed as a function of the
balance pattern. 
\cite{QL06} derived the
relationship between DD and the balance pattern for
symmetric designs. Furthermore, they calculated a new
lower bound of discrepancies
by using these connections.
In this section, we give similar results on the QQD criterion.
To our knowledge, this is the first time
that the balance pattern of an asymmetric design is studied,
and the extension from the symmetric situation to the
asymmetric one is not trivial.

Consider an asymmetric design $\D\in \cu(n,s_1^ps_2^q)$,
where the first $p$ factors
are qualitative and the last $q$ factors are quantitative.
For each $k$ columns of $\D$, $\d_{l_1},\dots,\d_{l_k}$
such that $1\leq l_1<\dots< l_{k_1}\leq p<
l_{k_1+1}<\dots< l_{k}\leq p+q$,
modify the definition of $B_{l_1,\dots,l_k}$ in \cite{FLW03}
by \bea\label{B_l1_lk}
B_{l_1,\dots,l_k}=\sum_{\Delta}
\left(n_{a_1,\dots,a_k}^{(l_1,\dots,l_k)}-
\frac{n}{s_1^{k_1}s_2^{k_2}}\right)^2,
\eea
to make it applicable to asymmetric situation,
where $\Delta=\{(a_1,\dots,a_k)|1\leq a_1,\dots,a_{k_1}\leq s_1,
1\leq a_{k_1+1},\dots, a_k \leq s_2\}, k=k_1+k_2$,  and
$n_{a_1,\dots,a_k}^{(l_1,\dots,l_k)}$ is the number of rows
in which the column group $\{\d_{l_1},\dots,\d_{l_k}\}$
takes the level combination $\{a_1,\dots,a_k\}$.
The summation is taken across all the possible level combinations.
If $B_{l_1,\dots,l_k}=0$, the sub-design $\{\d_{l_1},\dots,\d_{l_k}\}$
is an orthogonal array with strength $k$.
Obviously, $B_{l_1,\dots,l_k}$ assesses the closeness to
the orthogonality of strength $k$ of the sub-design
formed by the columns $\d_{l_1},\dots,\d_{l_k}$.
Now, we define the \textit{balance pattern}  by
$B(\D)=(B_1(\D),\dots,B_{p+q}(\D))$,
where \bea\label{Bk}
B_k(\D) = \sum_{\Omega} B_{l_1,\dots,l_k}\Bigg/\dbinom{p+q}{k},
\eea
with $\Omega=\{(l_1,\dots,l_k)| 1\leq l_1\cdots \leq l_{k_1}\leq p
< l_{k_1+1}\cdots\leq l_k \leq p+q\}.$
The summation is taken over all the possible $k$ columns.
As a result, $B_k(\D)$ evaluates the nearness between
the orthogonality of strength $k$ and the design
$\D$, and $B_k(\D)=0$ implies that the design $\D$ is
an orthogonal array of strength $k$.

The $B(\D)$ is defined based on the columns of a design,
interestingly, the following lemma illustrates that
it has a close relationship with the rows of a design.

\begin{lemma}\label{le_B_k}
For a design $\D \in \cu(n,s_1^ps_2^q)$,
the entries of $B(\D)$ could be reshaped as
\begin{align}\label{b_k}
B_k(\D)=\sum_{i,j=1}^{n}\sum_{\Omega}
\delta_{ij}^{(l_1,\dots,l_k)}\Bigg/\dbinom{p+q}{k}
-\sum_{\Omega}\frac{n^2}{s_1^{k_1}s_2^{k_2}}\Bigg/\dbinom{p+q}{k},
\text{ for }k = 1, \dots, p+q,
\end{align} where $\Omega=\{(l_1,\dots,l_k)| 1\leq l_1\cdots \leq l_{k_1}\leq p < l_{k_1+1}\cdots\leq l_k \leq p+q\},$
$\delta_{ij}^{(l_1,\dots,l_k)}$ equals $1$ if
$(x_{il_1},\dots,x_{il_k})=(x_{jl_1},\dots,x_{jl_k})$, and $0$ otherwise.
\end{lemma}

The proof of this lemma is very similar to
that of Lemma 3.1 in \cite{FLW03}.
The main difference is that $\D \in \cu(n,s_1^ps_2^q)$
is an asymmetric design when
$s_1 \neq s_2$. Thus we omit it.
Note that, the QQD in (\ref{qd_ex}) is defined by rows.
Now we can obtain another form of expression of QQD,
which relies on the balance pattern.

\begin{theorem}\label{qd_bp}
 If $\D=(\D_1,\D_2)$ with $\D_1 \in\cu(n,s^p)$
and $\D_2 \in\cu(n,2^q)$, then
\bea\label{qd_balance}
\text{QQD}^2(\D) = -\left(\frac{5s+1}{4s}\right)^p
\left(\frac{4}{3}\right)^q+\left(\frac{5s+1}{4s}\right)^p
\left(\frac{11}{8}\right)^q
+\frac{1}{n^2}\left(\frac{5}{4}\right)^{p+q}
\sum_{k=1}^{p+q}\left(\frac{1}{5}\right)^k
\dbinom{p+q}{k}B_k(\D).
\eea
\end{theorem}

Theorem \ref{qd_bp} reveals
the QQD is a function of $B_k(\D)$,
based on which we can derive a
new lower bound for the QQD.

\begin{theorem}\label{th_lb_bp}
 If $\D=(\D_1,\D_2)$ with $\D_1 \in\cu(n,s^p)$
and $\D_2 \in\cu(n,2^q)$,
then $QQD^2(\D)\geq \text{LB}_2$, where
\begin{align*}
LB_{2} =& -\left(\frac{5s+1}{4s}\right)^p \left(\frac{4}{3}\right)^q
+\left(\frac{5s+1}{4s}\right)^p \left(\frac{11}{8}\right)^q
+ \frac{1}{n^2}\left(\frac{5}{4}\right)^{p+q}
\sum_{k=1}^{p+q}\left(\frac{1}{5}\right)^k \\
&\times \sum_{k_1+k_2=k}
\dbinom{p}{k_1}\dbinom{q}{k_2}r_{n,k_1,k_2,s,2}
\left(1-\frac{r_{n,k_1,k_2,s,2}}{s^{k_1}2^{k_2}}\right),
\end{align*}
with $r_{n,k_1,k_2,s,2}$ being the residual of
$n~(\text{mod} ~ s^{k_1}2^{k_2})$.
\end{theorem}

For those designs $\D\in \cu(n,s^p2^q)$, two
lower bounds, $LB_1$ and $LB_2$, are applicable,
then define $LB=\max\{LB_1,LB_2\}$.

\begin{remark}
Let $p=0$ or $q=0$,
then the results related to QQD in
Theorems \ref{qd_bp} and \ref{th_lb_bp}
could degenerate into the corresponding results
in terms of  WD in \cite{FLW03} or
DD in \cite{QL06}, respectively.
\end{remark}

\section{Illustrative examples}\label{egs}

In this section, we give some examples to demonstrate
the reasonability of QQD for measuring the uniformity of designs
associated with both qualitative and quantitative factors.
Particularly, it turns out that the proposed QQD could distinguish
different MCDs generated by the same construction procedure.
In addition, it is shown that the design generated by
column juxtaposition of two separately uniform sub-designs
may not be uniform over the design space.

Many construction methods of MCDs were proposed\revA{.}
For example, \cite{HLS17} gave several
effective construction methods for MCDs with some
appealing space-filling properties for the quantitative factors
from the view of stratification on grids.
The generated MCDs are not unique because of
some structures and randomness.
Even if the MCDs share the same stratifications
on grids for the sub-designs of quantitative factors,
they may have distinct uniformity.
This declaration is interpreted by three MCDs visually.

\begin{example}\label{eg:MCD}
Construct \revA{an} MCD $\D_M = (\D_{M1},\D_{M2})$ via
Construction 3 in \cite{HLS17}
under the parameter settings: $s = 2, u = 4, k = 2$.
Then $\D_{M1}$ is an OA($16, 4, 2, 2$) and
$\D_{M2}$ is an LHD($16, 6$)
based on Theorem 5 in \cite{HLS17}.
Assume there are only one
qualitative and two quantitative factors in the design
system.

Now we consider selecting three columns from the large MCD,
$\D_{M}$, to form a small MCD, $\D = (\D_1,\D_2)$ with
one qualitative factor and two quantitative factors.
The $\D_1$ with one qualitative factor should be chosen from
$\D_{M1}$, 
such as its $1$st column. 
For the quantitative part $\D_2$, which should be a sub-design of
$\D_{M2}$, there are $\dbinom{6}{2}=15$ different choices.
For explanation, let $\D_{M2}=(\d_1,\dots,\d_6)$ and 
compare the following three situations,
$\D_2^{(1)}=(\d_2,\d_3)$, $\D_2^{(2)}=(\d_2,\d_4)$, and
$\D_2^{(3)}=(\d_3,\d_4)$, where
$\D_1 = (0 , 1 , 0 , 1 , 0 , 1 , 0 , 1 , 0 , 1 , 0 , 1 , 0 , 1 , 0 , 1)^T$,
$\d_2 = (0 , 1 , 15 , 14 , 8 , 9 , 6 , 7 , 11 , 10 , 5 , 4 , 2 , 3 , 13 , 12)^T$,
$\d_3 = (1 , 0 , 14 , 15 , 3 , 2 , 12 , 13 , 6 , 7 , 9 , 8 , 5 , 4 , 10 , 11)^T$
and
$\d_4 = (0 , 14 , 15 , 1 , 10 , 5 , 4 , 11 , 12 , 3 , 2 , 13 , 6 , 9 , 8 , 7)^T$.
Denote the resulting three small MCDs by
$\D^{(1)} = (\D_1, \D_2^{(1)})$, $\D^{(2)} = (\D_1, \D_2^{(2)})$ and
$\D^{(3)} = (\D_1, \D_2^{(3)})$. It can be easily calculated that
$\text{QQD}^2(\D^{(1)}) = 0.0066$, $\text{QQD}^2(\D^{(2)}) = 0.0063$
and $\text{QQD}^2(\D^{(3)}) = 0.0060$.
Thus, under the QQD criterion,
$\D^{(3)}$ is the most uniform design among the
three designs,
followed by $\D^{(2)}$, and $\D^{(1)}$ is the worst one.

Figure \ref{Fig:MCD} displays the scatter plots of the three
sub-designs $\D_2^{(1)}$, $\D_2^{(2)}$ and $\D_2^{(3)}$
for the quantitative factors. According to Theorem 5 in \cite{HLS17},
$\D_2^{(1)}$ achieves stratification on $2 \times 2$ grids,
and both $\D_2^{(2)}$ and $\D_2^{(3)}$ can achieve
stratification on $8 \times 2$ grids and $2 \times 8$ grids.
Consequently, $\D_2^{(1)}$ is the worst among the three sub-designs,
and $\D_2^{(2)}$ and $\D_2^{(3)}$ possess the same performance
in terms of space-filling property on grids.
We further consider the points of $\D_2^{(i)},i=1,2,3$,
corresponding to level $0$ and $1$ of $\D_1$, respectively,
which form a small LHD since $\D^{(i)}$ are all MCDs.
Additionally, the points of $\D_2^{(1)}$ and $\D_2^{(3)}$
corresponding to any level of $\D_1$ can achieve stratification on
$2 \times 2$ grids, but there are clustered points in $\D_2^{(1)}$.
From Figure \ref{Fig:MCD}, we can conclude that
$\D^{(3)}$ is the most outstanding, followed by $\D^{(2)}$
and $\D^{(1)}$ is the worst, which is coincident with the results
in comparison under the QQD criterion.
\begin{figure}[h!]
\centering
\includegraphics[width=0.96\textwidth]{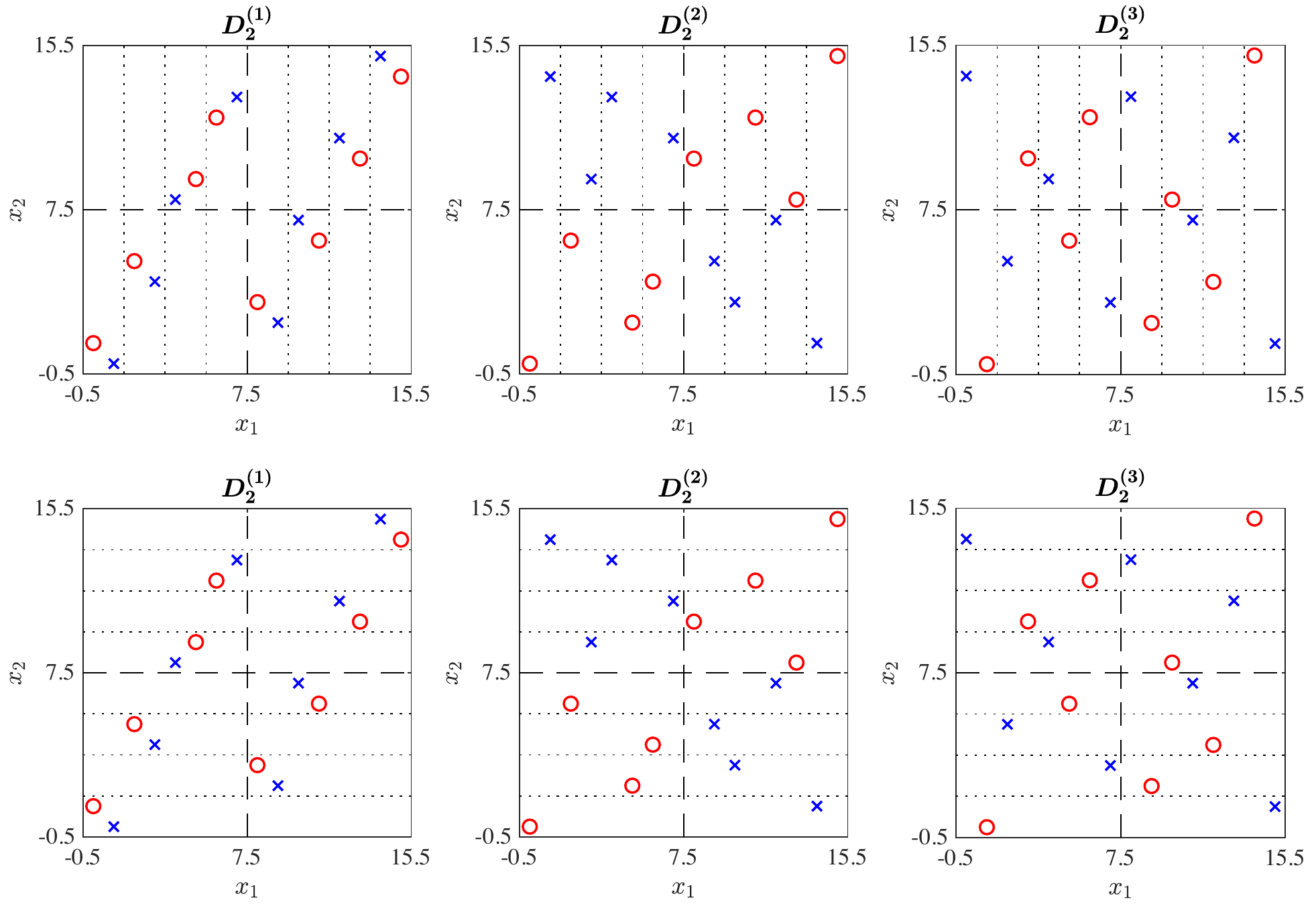}
\caption{The scatter plots for the quantitative sub-designs,
$\D_2^{(1)}$, $\D_2^{(2)}$ and $\D_2^{(3)}$ in Example
\ref{eg:MCD}. The points of $\D_2^{(j)} ( j = 1, 2 , 3)$
corresponding to level 0 and 1 of $\D_1$ are represented by
``\revA{$\circ$}" and ``\revA{\footnotesize{$\times$}}", respectively.
The panels in the first row shows
stratification on $8\times 2$ grids and the panels in the
second row exhibits stratification on $2\times 8$ grids.
}\label{Fig:MCD}
\end{figure}
\end{example}

Example \ref{eg:MCD} reflects that the QQD can distinguish
the MCDs constructed by the same construction method
even though they have the same stratifications
on grids for the quantitative factors,
like the $\D^{(2)}$ and $\D^{(3)}$.
Originally, although the proposed QQD is motivated by
comparing different MCDs, it is not limited to MCDs,
and can be used for quantifying the uniformity
of any design with both types of factors.

Since $\D$ consists of two sub-designs, $\D_1$ and
$\D_2$, one may consider searching uniform designs
$\D_1^*$ and $\D_2^*$ for two sub-designs, separately,
then juxtaposing them by column to obtain
$\D^*=(\D_1^*,\D_2^*)$. However, we find some
unsatisfactory results about this method.

\begin{example}\label{eg:2+2}
If $n=16,p=2,q=2,s_1=s_2=2,s_3=s_4=4$, consider
column juxtaposition of two separately uniform designs.

Let $\D_1$ and $\D_2^{(1)}$ be the full factorial designs as below.
Permute the rows in $\D_2^{(1)}$ and obtain another $4$-level
full factorial design $\D_2^{(2)}$, where
$$\D_1 = \left(\begin{array}{cccccccccccccccc}
0 & 0 & 0 & 0 & 0 & 0 & 0 & 0 & 1 & 1 & 1 & 1 & 1 & 1 & 1 & 1 \\
0 & 0 & 0 & 0 & 1 & 1 & 1 & 1 & 1 & 1 & 1 & 1 & 0 & 0 & 0 & 0
\end{array}\right)^T,$$
$$\D_2^{(1)} = \left(\begin{array}{cccccccccccccccc}
0 & 0 & 0 & 0 & 1 & 1 & 1 & 1 & 2 & 2 & 2 & 2 & 3 & 3 & 3 & 3 \\
0 & 1 & 2 & 3 & 0 & 1 & 2 & 3 & 0 & 1 & 2 & 3 & 0 & 1 & 2 & 3
\end{array} \right)^T,$$ and
$$\D_2^{(2)} = \left(\begin{array}{cccccccccccccccc}
0 & 0 & 1 & 1 & 2 & 2 & 3 & 3 & 0 & 0 & 1 & 1 & 2 & 2 & 3 & 3 \\
0 & 2 & 1 & 3 & 0 & 2 & 1 & 3 & 1 & 3 & 0 & 2 & 1 & 3 & 0 & 2
\end{array} \right)^T.$$
Obviously, $\D_1$, $\D_2^{(1)}$ and $\D_2^{(2)}$ are
uniform designs, respectively. Let
$\D^{(1)} = (\D_1, \D_2^{(1)})$ and $\D^{(2)} = (\D_1, \D_2^{(2)})$,
that is, both resulting designs $\D^{(1)}$ and $\D^{(2)}$
consist of two separately uniform designs.
Figure \ref{Fig:2+2} shows the scatter points.
Scatter points for the quantitative factors of both
$\D^{(1)}$ and $\D^{(2)}$ have the same uniformity.
However, we can see that the points of $\D_2^{(2)}$
corresponding to level $0$ and $1$ of each factor of $\D_1$
have better uniformity than that of $\D_2^{(1)}$, namely,
the points ``\revA{$\circ$}'' and
``\revA{\footnotesize{$\times$}}" in $\D_2^{(2)}$
is more uniform than the points ``\revA{$\circ$}'' and
``\revA{\footnotesize{$\times$}}" in $\D_2^{(1)}$, respectively.
For instance, with respect to the level $0$ of the first factor of
$\D_1$, the corresponding points of $\D_2^{(1)}$, i.e.,
points ``\revA{$\circ$}'' in the left top panel, only explore
half of the region of $x_1$, while
the points of $\D_2^{(2)}$ corresponding to the level $0$
of the first factor of $\D_1$ i.e.,
points ``\revA{$\circ$}'' in the right top panel, can probe
the whole region of $x_1$. Therefore, we can assert that
$\D^{(2)}$ looks more uniformly than $\D^{(1)}$ from the
scatter points.

On the other hand, from Theorem \ref{th_qd}, we have
$\text{QQD}^2(\D^{(1)}) = 0.0822$ and
$\text{QQD}^2(\D^{(2)}) = 0.0545$, which implies
$\D^{(2)}$ is more uniform than $\D^{(1)}$ under the QQD criterion, agreeing with our intuition.

\begin{figure}[h!]
\centering
\includegraphics[width=0.7\textwidth]{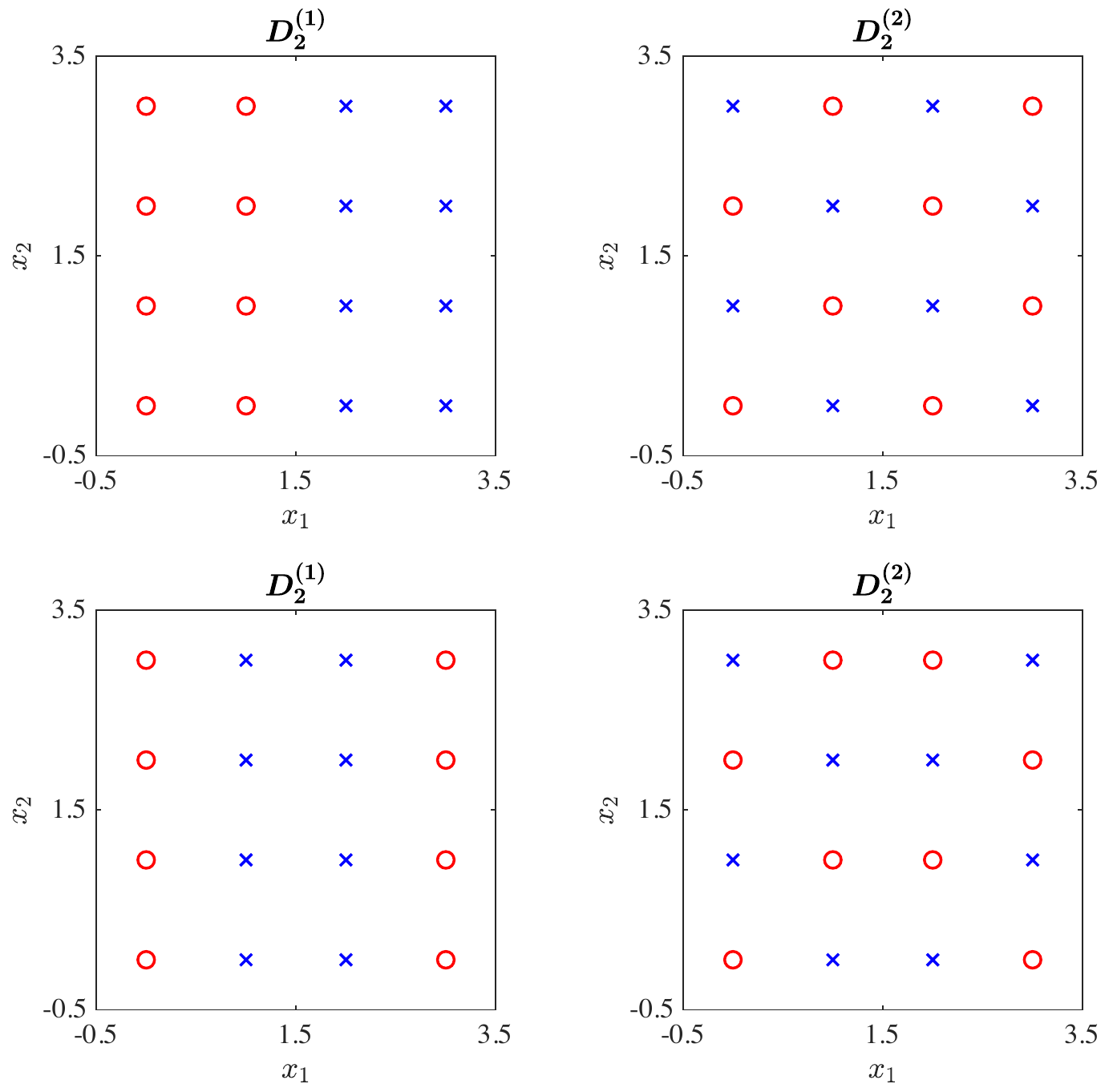}
\caption{The scatter plots for the quantitative factors
$\D_2^{(1)}$ and $\D_2^{(2)}$ in Example \ref{eg:2+2}.
The points of $\D_2^{(j)}, j = 1, 2$, corresponding to
level $0$ and $1$ of $\D_1$ are represented by
``\revA{$\circ$}" and ``\revA{\footnotesize{$\times$}}",
respectively. The top and bottom panels show plots with
respect to the first and the second factor of $\D_1$,
respectively. }\label{Fig:2+2}
\end{figure}
\end{example}

This example not only shows
that the proposed QQD can discriminate different designs,
but also indicates that column juxtaposition
of two separately uniform designs does not work for
constructing the uniform designs with both types of
factors since the combination mode of the
qualitative and quantitative factors could affect the
uniformity of the joint designs.

For comparison, a naive criterion for designs with both qualitative and quantitative factors, the sum of the space-filling criterion of the quantitative factors over each level of each qualitative factor, may be considered.
Such \revA{a} criterion only measures the space-filling properties of all the quantitative factors and the space-filling properties between each qualitative factor and all the quantitative factors.
Consider the sum of the WD values of the quantitative factors over each level of each qualitative factor as one naive criterion, denoted by SWD. Then SWD($\D^{(2)}$) sums the WD values of the four 8-run sub-designs, represented by ``\revA{$\circ$}" in the right top and right bottom panel, and  ``\revA{\footnotesize{$\times$}}" in the right top and the right bottom panel, respectively. It can be easily calculated that SWD($\D^{(2)}$) = $1.1055$. Now introduce another design for comparison, $\widetilde{\D}^{(2)}=(\widetilde{\D}_1,\D_2^{(2)})$, where $\widetilde{\D}_1 = \left(\zero_{2\times8},\one_{2\times8}
\right)^T$. Obviously, $\D^{(2)}$ has better space-filling property
than $\widetilde{\D}^{(2)}$ since the two columns for the qualitative factors in $\widetilde{\D}^{(2)}$ remain the same. By computation, SWD($\widetilde{\D}^{(2)}$) = $1.0999 < $ SWD($\D^{(2)}$) implies that $\widetilde{\D}^{(2)}$ is better than $\D^{(2)}$ under SWD. This result is contrary to the objective truth.
In terms of QQD criterion, it could be calculated that $\text{QQD}^2(\widetilde{\D}^{(2)}) = 0.0813 > \text{QQD}^2(\D^{(2)})$. The comparison result under QQD is still agrees with the facts. It confirms that the advantage of QQD over the naive criterion in the projection space-filling properties on all subspaces are important.

Next, we give two designs, whose QQD
values achieve the lower bounds obtained in this paper.
This further indicates the lower bounds we derived are tight.

\begin{example}\label{ex_lb1}
Consider the design $\D=(\D_1,\D_2)\in \cu(4,4\times2^2)$ with
$\D_1=(0,1,2,3)^T$ and
$\D_2 = \big((0,1,0,1)^T, ~(1,0,0,1)^T\big)$.
It reaches the lower bound $LB_2$ with $\text{QQD}^2(\D)=
LB_2= 0.1706$.
\end{example}

\begin{example}\label{ex_lb2}
Consider the design $\D$ listed in Table \ref{ta_ex_lb2},
which can \revA{reach} the lower bound $LB_1$, and
$\text{QQD}^2(\D)=LB_1= 17.0235$.
\begin{table}[!h]
\caption{$\D=(\D_1,\D_2)\in \cu(8,2^7\times4^7)$ in
Example \ref{ex_lb2}}\label{ta_ex_lb2}  
\bc {\small \tabcolsep=8pt \bt {c|c} \hline $\D_1$ & $\D_2$
\\\hline  $ \ba{ccccccc}

0&	0&	0&	0&	1&	0&	0\\
0&	0&	0&	1&	0&	1&	1\\
0&	1&	1&	1&	0&	0&	0\\
1&	1&	0&	1&	1&	0&	1\\
0&	1&	1&	0&	1&	1&	1\\
1&	0&	1&	1&	1&	1&	0\\
1&	0&	1&	0&	0&	0&	1\\
1&	1&	0&	0&	0&	1&	0

 \ea  $ &  $\ba{ccccccc}
1&	1&	2&	3&	1&	2&	3\\
3&	2&	0&	0&	1&	1&	2\\
3&	3&	1&	2&	2&	3&	3\\
2&	2&	2&	2&	0&	0&	1\\
2&	3&	3&	0&	3&	2&	0\\
1&	0&	3&	1&	2&	0&	2\\
0&	1&	0&	1&	0&	3&	0\\
0&	0&	1&	3&	3&	1&	1
 \ea $  \\
 \hline
 \et} \ec
 \end{table}
\end{example}

The designs presented in Examples \ref{ex_lb1} and \ref{ex_lb2}
show that the
lower bounds $LB_1$ and $LB_2$ obtained in this paper
are tight 
and achievable.
\revA{
Furthermore, we give the following example to illustrate the use and effectiveness of the QQD in comparing response surface designs.

\begin{example}\label{RSD}
Consider a response surface design with one qualitative ($z$) and two quantitative factors ($x_1$ and $x_2$). Table \ref{ta_ex_RSD} lists the four designs, among which the quantitative factors are fixed.
The sub-design composed by $x_1$ and $x_2$ is a central composite design (CCD), which consists of a $2^2$ factorial design, two center runs and four axial runs. The difference among the four designs exists in the assignment for the qualitative factor. The $z^{(1)}, \dots, z^{(4)}$ in Table \ref{ta_ex_RSD} are respectively represented as the  qualitative factor $z$ in the four designs.
Denote the four designs by $\D^{(1)} = (x_1,x_2,z^{(1)})$, $\D^{(2)} = (x_1,x_2,z^{(2)})$, $\D^{(3)} = (x_1,x_2,z^{(3)})$, and $\D^{(4)} = (x_1,x_2,z^{(4)})$, respectively.
\cite{WD98} compared the efficiency of different designs by using the determinant criterion (D-criterion).
They obtained that the first design $\D^{(1)}$ is the optimal because of its maximum D-criterion value.
It is noted that the D-criterion heavily relies on the pre-determined model.
A small misspecification may seriously affect the performance of the chosen design. Then, the space-filling property may be considered.
Next, we will choose a new design by the QQD criterion.

\begin{table}[!h]
\renewcommand\thetable{2}
\renewcommand\arraystretch{1.3}
\caption{The values of $x_1$, $x_2$, and $z$ in the four designs in
Example 6}\label{ta_ex_RSD}
\bc {\small \tabcolsep=8pt
\begin{tabular}{c|rr|rrrr}\hline
Run & ~$x_1$ & $x_2$ & $z^{(1)}$ & $z^{(2)}$ &  $z^{(3)}$  & $z^{(4)}$ \\
\hline
1  & 1 & 1 & $-1$ & $-1$ & $-1$ & $-1$~ \\
2  & 1 & $-1$ & $-1$ & 1  & $-1$ & 1~  \\
3  & $-1$  & 1  & $-1$  & 1  & 1  & 1~  \\
4  & $-1$ & $-1$  & 1   & 1   & 1   & $-1$~ \\ \hdashline
5  & 0 & 0 & 1   & 1   & 1   & 1~  \\
6  & 0 & 0 & $-1$ & $-1$ & $-1$ & $-1$~ \\ \hdashline
7  &$\sqrt{2}$ & 0  & 1   & $-1$ & $-1$ & $1$~ \\
8  &$-\sqrt{2}$ & 0 & $-1$ & $-1$ & $-1$ & $-1$~  \\
9  & 0 & $\sqrt{2}$ & 1  & $-1$ & $-1$ & $-1$~ \\
10  & 0 & $-\sqrt{2}$  & $-1$ & $-1$ & 1  & 1~  \\ \hline
\end{tabular}} \ec
\end{table}

According to the structure of the CCD, the sub-design composed by the two quantitative factors consists of three portions, factorial, center, and axial portion, see the three portions of $(x_1,x_2)$ listed in Table \ref{ta_ex_RSD}.
Corresponding to the factorial portion of CCD,
denote the sub-designs composed by the first four runs of the four designs by $\D_1^{(1)}$, $\D_1^{(2)}, \D_1^{(3)}$ and $\D_1^{(4)}$, respectively.
When the space-filling property of the factorial portion is considered, we can use the QQD to assess the uniformity of the four sub-designs $\D_1^{(1)}, \dots, \D_1^{(4)}$.
After transforming the level $-1$ to $1/4$ and $1$ to $3/4$ to make the values of $x_1$ and $x_2$ falling into [0,1],
the QQD values could be calculated according to the expression of QQD. Table \ref{ta_ex_RSD_value} lists the QQD values of the four sub-designs. It can be seen that $\D_1^{(4)}$ owes the minimum QQD value. Moreover, QQD$(\D_1^{(4)})$ reaches the lower bound of QQD for designs in $\cu(4,2^3)$. These results imply that $\D_1^{(4)}$ is the best among the four sub-designs under QQD, and the four sub-designs are ranked as $\D_1^{(4)}, \D_1^{(3)}, \D_1^{(1)}, \D_1^{(2)}$.
Next, consider the space-filling property of the four full designs, $\D^{(1)}, \dots, \D^{(4)}$. At this time, all the four designs are not U-type.
Transform each level of $x_1$ and $x_2$ into $[0,1]$ by
$\big(x+\sqrt{2}\big)/\big(2\sqrt{2}\big)$, $x\in\big\{-\sqrt{2}, -1, 0, 1, \sqrt{2}\big\}$, then, we can compute the QQD values of $\D^{(1)}, \dots, \D^{(4)}$. Based on the QQD values listed in Table \ref{ta_ex_RSD_value}, $\D^{(4)}$ is the best design.
Combining the comparison result of the four sub-designs, the fourth design is suggested when the space-filling property is considered.

Additionally, the extended MaxPro criterion proposed in \cite{J19} could be considered here to assess the space-filling property. The corresponding criterion values of the four sub-designs $\D_1^{(1)}, \dots, \D_1^{(4)}$, and the four  designs $\D^{(1)}, \dots, \D^{(4)}$, are also listed in Table \ref{ta_ex_RSD_value}. It can be seen that, under the  extended MaxPro criterion, the order of the four sub-designs is accordant with that under the QQD, and the best design among the four designs $\D^{(1)}, \dots, \D^{(4)}$ is the same with that under the QQD. 
It is noted that there is a small difference of the order of $\D^{(2)}$ and $\D^{(3)}$ under the two criteria. It means that the two criteria have a litter bit difference,
and which does not affect the  suggested design.

\begin{table}[!t]
\caption{The criterion values of designs under two space-filling criteria}\label{ta_ex_RSD_value}
\bc {\small \tabcolsep=7pt
\begin{tabular}{c|cccc|cccc}\hline
Criterion & $\D_1^{(1)}$ & $\D_1^{(2)}$ & $\D_1^{(3)}$ & $\D_1^{(4)}$
& $\D^{(1)}$ & $\D^{(2)}$ & $\D^{(3)}$ & $\D^{(4)}$
\\ \hline
Squared QQD       &
0.2255 & 0.2255 & 0.1766 & 0.1571 &
0.0763 & 0.0795 & 0.0792 & 0.0653 \\
Extended MaxPro &
1.8821 & 1.8821 & 1.8246 & 1.3606 &
3.5810 & 3.6441 & 3.7086 & 3.5587 \\ \hline
\end{tabular}} \ec
\end{table}

\end{example}
}

\section{Conclusion}
For measuring the space-filling property of designs with both
qualitative and quantitative factors, in this paper,
we propose a new uniformity criterion,
the qualitative-quantitative discrepancy, and give its explicit expression.
From several intuitive examples, we show the reasonability
and effectiveness of the proposed QQD.
In fact, the \revA{newly proposed} discrepancy can measure the uniformity of any design containing two types of factors apart from that of MCDs \revA{and response surface designs when the space-filling properties are considered.}
In addition, according to the closed form of the QQD and
the connections between QQD and the balance pattern,
two tight lower bounds of QQD are obtained
by strict mathematic deductions, which can
be regarded as a benchmark for identifying the
uniform designs.
\revA{When the QQD is applied to stochastic optimization algorithms,
as similar as the CD, WD and MD,
an iteration formula could be derived to greatly simplified the calculation of the QQD value for a new design in each iteration.}
In the future, we can
use some stochastic optimization methods or some systematic
construction methods to search uniform designs
under QQD criterion.

\section*{Acknowledgements}
This work was supported by National Natural Science Foundation of China (11871288) and Natural Science Foundation of Tianjin (19JCZDJC31100). The first two authors contributed equally to this work.

\section*{Appendix}
\noindent \textbf{\textit{Proof of Theorem \ref{th_qd}}}.
Note that $\chi=\chi_1\times\dots\times\chi_{p+q}$, where
$\chi_k=\{0,1,\dots,\s_k-1\}$, for $k=1,2,\dots,p$ and
$\chi_k=[0,1]$, for $k=p+1,\dots,p+q$.
For the first term of (\ref{discre}), we have
\bea\label{term1}
\int_{\chi_k^{2}}\ck_k(t_k,z_k)\mathrm{d}F(t_k) \mathrm{d}F(z_k) =
\begin{cases}
 \frac{5s_k+1}{4s_k}, & \text{for $k=1,\dots,p$},\\
\frac{4}{3}, & \text{for $k=p+1,\dots,p+q.$}
\end{cases}
\eea
Hence, $\int_{\chi^{2}}\ck(\t,\z)\mathrm{d}F(\t) \mathrm{d}F(\z)
=\prod_{k=1}^{p}\left(\frac{5s_k+1}{4s_k}\right)\left(\frac43\right)^q$.
For the second term,
\bea\label{term2}
\int_{\chi_k}\ck_k(t_k,x_{ik})\mathrm{d}F(t_k)=
\begin{cases}
\frac{5s_k+1}{4s_k}, & \text{for $k=1,\dots,p$},\\
\frac{4}{3}, & \text{for $k=p+1,\dots,p+q.$}
\end{cases}
\eea
Thus, $\frac{2}{n}\sum_{i=1}^{n}
\int_{\chi}\ck(\t,\x_i)\mathrm{d}F(\t)=2\prod_{k=1}^{p}
\left(\frac{5s_k+1}{4s_k}\right)\left(\frac{4}{3}\right)^q$.
With regard to the last term,
\bea\label{term3}
\sum_{i,j=1}^{n}\ck(\x_i,\x_j)=\sum_{i,j=1}^{n}
\left(\frac{3}{2}\right)^{\delta_{ij}(\D_1)}\left(\frac{5}{4}\right)^{p-\delta_{ij}(\D_1)}
\cdot\prod_{k=p+1}^{p+q}\left(
\frac{3}{2}-|~x_{ik}-x_{jk}~|+|~x_{ik}-x_{jk}~|^2\right).
\eea
Substituting the Equations (\ref{term1}) - (\ref{term3})
into (\ref{discre}), we complete the proof.
$\hfill\blacksquare$

For proving Corollary \ref{co_full}, we give a lemma,  which is  the extension of related results in \cite{ZFN12} and its proof is  straightforward.
\begin{lemma}\label{le_property}
The matrices $\A$ and $\A_k$ in Lemma \ref{quadratic_form} have the following properties,
\begin{align*}
\begin{cases}
~\A_k\one_{s_k}=\left(\frac{3}{2}+\frac{5}{4}(s_k-1)\right)\one_{s_k}, \\
~\A_k^{-1}\one_{s_k}=\left(\frac{3}{2}+\frac{5}{4}(s_k-1)\right)^{-1}\one_{s_k},
\end{cases}\text{for}~ k=1,\dots,p,~~ \end{align*}\begin{align*}
\begin{cases}
~\A_k\one_{s_k}= \left(\frac{4s_k}{3}+\frac{1}{6s_k}\right)\one_{s_k}, \\
~\A_k^{-1}\one_{s_k}= \left(\frac{4s_k}{3}+\frac{1}{6s_k}
\right)^{-1}\one_{s_k},
\end{cases} \text{for}~ k=p+1,\dots,p+q,
 \end{align*}
  and
 \begin{align*}
 \begin{cases}
 \A\one_{N}=\prod_{k=1}^p\left(\frac{3}{2}+\frac{5}{4}(s_k-1)\right)\prod_{k=p+1}^{p+q}
 \left(\frac{4s_k}{3}+\frac{1}{6s_k}\right)\one_{N},\\
  \A^{-1}\one_{N} =\prod_{k=1}^p\left(\frac{3}{2}+\frac{5}{4}(s_k-1)\right)^{-1}\prod_{k=p+1}^{p+q}
  \left(\frac{4s_k}{3}+\frac{1}{6s_k}\right)^{-1}\one_{N}.
 \end{cases}
 \end{align*}
\end{lemma}

\vspace{2mm} \noindent \textbf{\textit{Proof of Corollary \ref{co_full}}}. From Lemma \ref{quadratic_form} and Theorem \ref{th_uni}, we have
\begin{align}\label{proof_co_1}
\text{QQD}^2(\D^*) &= -\prod_{k=1}^{p}\left(\frac{5s_k+1}{4s_k}\right)
\left(\frac{4}{3}\right)^q+\frac{1}{n^2}\left(\frac{n}{N}\one_N\right)^T\A\left(\frac{n}{N}\one_N\right)\nonumber
\\
&=-\prod_{k=1}^{p}\left(\frac{5s_k+1}{4s_k}\right)
\left(\frac{4}{3}\right)^q+ \frac{1}{N^2}\one_N^T\A\one_N
\end{align}
Note that $N= \prod_{k = 1}^{p+q}s_k$ and by Lemma \ref{le_property},
 \begin{align}\label{proof_co_2}
 \one_{N}^T\A\one_{N}&=\prod_{k=1}^p\left(\frac{3}{2}+\frac{5}{4}(s_k-1)\right)\prod_{k=p+1}^{p+q}
 \left(\frac{4s_k}{3}+\frac{1}{6s_k}\right)\one_{N}^T\one_{N}
 \nonumber\\
 &=N\prod_{k=1}^p\left(\frac{3}{2}+\frac{5}{4}(s_k-1)\right)\prod_{k=p+1}^{p+q}
 \left(\frac{4s_k}{3}+\frac{1}{6s_k}\right)  \nonumber\\
 &=N^2\prod_{k=1}^p\left(\frac{5s_k+1}{4s_k}\right)\prod_{k=p+1}^{p+q}
 \left(\frac{4}{3}+\frac{1}{6s^2_k}\right)
\end{align}
Combining (\ref{proof_co_1}) and (\ref{proof_co_2}), this proof can be finished.
$\hfill\blacksquare$

\vspace{2mm} \noindent \textbf{\textit{Proof of Theorem \ref{th_lbqd}}}.
Note that $\left(\frac{6}{5}\right)^{\delta_{ij}(\D_1)}=\prod_{k=1}^p
\left(\frac{6}{5}\right)^{\delta_{ij}^k}$,
where $\delta_{ij}^k=\delta_{x_{ik}x_{jk}}$.
Let $\alpha_{ij}^k=|~x_{ik}-x_{jk}~|-|~x_{ik}-x_{jk}~|^2$,
for $k=p+1,\dots,p+q$.
The formula (\ref{qd_ex}) can be reshaped as
\begin{align*}
\text{QQD}^2(\D) = C+ \frac{1}{n}\left(\frac{3}{2}\right)^{p+q}
+\frac{1}{n^2}\left(\frac{5}{4}\right)^p\sum_{1\leq i\neq j\leq n}
\prod_{k=1}^p\left(\frac{6}{5}\right)^{\delta_{ij}^k}
\prod_{k=p+1}^{p+q}\left(\frac{3}{2}-\alpha_{ij}^k\right).
\end{align*}
From the above equation,
it's clear that the value of QQD is a function of
$\prod_{k=1}^p\left(\frac{6}{5}\right)^{\delta_{ij}^k}\prod_{k=p+1}^{p+q}
\left(\frac{3}{2}-\alpha_{ij}^k\right)$, $i,j=1,\dots,n$, $i\neq j$, as
the first two terms are constant for the given parameters.
For the U-type design $\D_1\in \cu(n,s_1\cdots s_p)$,
there are $\frac{n^2(s_k-1)}{s_k}$ number of
$\delta_{ij}^k$ being 0, and $\frac{n(n-s_k)}{s_k}$
number of $\delta_{ij}^k$
being 1, $i,j=1,\dots,n,i\neq j$, for any $k=1,\dots,p$.
For $\D_2\in \cu(n,s_{p+1} \cdots  s_{p+q})$,
Table \ref{ta_2} shows the distribution of
$\alpha_{ij}^k, i,j=1,\dots,n,i\neq j$, for any $k=p+1,\dots,p+q$.
Under the given parameters $(n,s_1,\dots,s_{p+q})$ of a design,
minimizing $\text{QQD}^2(\D)$ is equivalent to minimizing
$\sum_{1\leq i\neq j\leq n}\prod_{k=1}^p
\left(\frac{6}{5}\right)^{\delta_{ij}^k}\prod_{k=p+1}^{p+q}
\left(\frac{3}{2}-\alpha_{ij}^k\right)$.
From the above discussion, we claim that
$\prod_{1\leq i\neq j\leq n}\prod_{k=1}^p
\left(\frac{6}{5}\right)^{\delta_{ij}^k}\prod_{k=p+1}^{p+q}
\left(\frac{3}{2}-\alpha_{ij}^k\right)$ is a constant for the given
$(n,s_1,\dots,s_{p+q})$, say $H$.
According to the geometric and arithmetic mean inequality,
when each $\prod_{k=1}^p\left(\frac{6}{5}\right)^{\delta_{ij}^k}
\prod_{k=p+1}^{p+q}\left(\frac{3}{2}-\alpha_{ij}^k\right)$
takes the same values for $1\leq i\neq j \leq n$, i.e., 
$H^{\frac{1}{n(n-1)}}$,
$\text{QQD}(\D)^2$ can \revA{reach} its minimum, which completes
the proof.

\begin{table}[!h] \linespread{1.5}
\caption{Distribution of $\alpha_{ij}^k,\text{for} ~k=p+1,\dots,p+q$}
\label{ta_2}
\bc { \tabcolsep=10pt \large
\bt{cc|cc}
\cline{1-2} \cline{3-4}
\multicolumn{2}{c|}{even $s_k$} & \multicolumn{2}{c}{odd $s_k$} \\\hline
the values of $\alpha_{ij}^k$ & number &
the values of $\alpha_{ij}^k$ & number\\\hline
0 & $\frac{n(n-s_k)}{s_k}$ & 0 & $\frac{n(n-s_k)}{s_k}$
\\
$\frac{2(2s_k-2)}{4s_k^2}$ & $\frac{2n^2}{s_k}$
& $\frac{2(2s_k-2)}{4s_k^2}$ & $\frac{2n^2}{s_k}$ \\
\vdots & \vdots & \vdots & \vdots\\
$\frac{(s_k-2)(s_k+2)}{4s_k^2}$ & $\frac{2n^2}{s_k}$
& $\frac{(s_k-3)(s_k+3)}{4s_k^2}$ & $\frac{2n^2}{s_k}$\\
$\frac{s_k^2}{4s_k^2}$ & $\frac{n^2}{s_k}$
& $\frac{(s_k-1)(s_k+1)}{4s_k^2}$ & $\frac{2n^2}{s_k}$
 \\\hline
\et} \ec
\end{table}
$\hfill\blacksquare$

\vspace{2mm} \noindent \textbf{\textit{Proof of Theorem \ref{qd_bp}}}. For
$\D_1 \in\cu(n,s^p)$
and $\D_2 \in\cu(n,2^q)$, from Theorem \ref{th_qd},
we have
\begin{align}\label{pr_th_1}
\text{QQD}^2(\D)=&C
+\frac{1}{n^2}\left(\frac{5}{4}\right)^p\sum_{i,j=1}^{n}
\left(\frac{6}{5}\right)^{\delta_{ij}(\D_1)}\cdot\prod_{k=p+1}^{p+q}
\left(\frac{3}{2}-|~x_{ik}-x_{jk}~|+|~x_{ik}-x_{jk}~|^2\right)\nonumber\\
=&C
+\frac{1}{n^2}\left(\frac{5}{4}\right)^p\sum_{i,j=1}^{n}
\prod_{k=1}^p\left(\frac{6}{5}\right)^{\delta_{ij}^k}
\cdot\prod_{k=p+1}^{p+q}
\left(\frac{3}{2}-|~x_{ik}-x_{jk}~|+|~x_{ik}-x_{jk}~|^2\right),
\end{align}
where $C=-\left(\frac{5s+1}{4s}\right)^p\left(\frac{4}{3}\right)^q$.
 Note that,
$\left(\frac{6}{5}\right)^{\delta_{ij}^k} = 1+\frac15\delta_{ij}^k$,
for $k = 1, \dots p$, and
$|~x_{ik}-x_{jk}~| = \frac{1}{2}(1-\delta_{ij}^k)$, for
$k = p+1, \dots p+q$. Thus, (\ref{pr_th_1}) can
be rewritten as
\begin{align}\label{pr_th_4}
\text{QQD}^2(\D)=-\left(\frac{5s+1}{4s}\right)^p\left(\frac{4}{3}\right)^q
+\frac{1}{n^2}\left(\frac{5}{4}\right)^p
\sum_{i,j=1}^{n}\prod_{k=1}^p\left(1+\frac{1}{5}\delta_{ij}^k\right)
\cdot\prod_{k=p+1}^{p+q}\left(\frac{5}{4}+\frac{1}{4}\delta_{ij}^k\right).
\end{align}
The value of (\ref{pr_th_4}) is dominated by the last
term for given parameters, because the first term is a constant.
Part of the last term
\begin{align}
&\sum_{i,j=1}^{n}\prod_{k=1}^p\left(1+\frac{1}{5}\delta_{ij}^k\right)
\cdot\prod_{k=p+1}^{p+q}\left(\frac{5}{4}+\frac{1}{4}\delta_{ij}^k\right)
\nonumber\\
&= \sum_{i,j=1}^{n}\left(\left(\frac{5}{4}\right)^q
+\sum_{k=1}^{p+q}\sum_{\Omega}\left(\frac{1}{5}\right)^{k_1}
\left(\frac{1}{4}\right)^{k_2}\left(\frac{5}{4}\right)^{q-k_2}
\delta_{ij}^{(l_1,\dots,l_k)} \right) \nonumber\\
&=  n^2\left(\frac{5}{4}\right)^q+\left(\frac{5}{4}\right)^q
\sum_{k=1}^{p+q}\left(\frac{1}{5}\right)^k\sum_{i,j=1}^{n}\sum_{\Omega}
\delta_{ij}^{(l_1,\dots,l_k)} \nonumber \\ 
&= n^2\left(\frac{5}{4}\right)^q+\left(\frac{5}{4}\right)^q
\sum_{k=1}^{p+q}\left(\frac{1}{5}\right)^k
\left(B(k)+\sum_{\Omega}\frac{n^2}{s^{k_1}2^{k_2}}\Bigg/
\dbinom{p+q}{k}\right)\dbinom{p+q}{k}\label{pr_th_5}\\
&= n^2\left(\frac{5}{4}\right)^q+ \left(\frac{5}{4}\right)^q
\sum_{k=1}^{p+q}\left(\frac{1}{5}\right)^k\dbinom{p+q}{k}B(k)
+ n^2\left(\frac{5}{4}\right)^q\sum_{k=1}^{p+q}
\left(\frac{1}{5}\right)^k\sum_{\Omega}\frac{1}{s^{k_1}2^{k_2}}.
\label{pr_th_6}
\end{align}
The Equation (\ref{pr_th_5}) holds because of (\ref{b_k})
in Lemma \ref{le_B_k}. Furthermore,
\begin{align}\label{pr_th_7}
\sum_{k=1}^{p+q}
\left(\frac{1}{5}\right)^k\sum_{\Omega}\frac{1}{s^{k_1}2^{k_2}}
=& \sum_{k=1}^{p+q}\sum_{\Omega}\frac{1}{(5s)^{k_1}10^{k_2}}
= \left(1+\frac{1}{5s}\right)^p \left(1+\frac{1}{10}\right)^q-1\nonumber\\
=& \left(\frac{5s+1}{5s}\right)^p\left(\frac{11}{10}\right)^q-1.
\end{align}
Substituting (\ref{pr_th_7}) into (\ref{pr_th_6}),
and substituting (\ref{pr_th_6}) into (\ref{pr_th_4}),
we obtain (\ref{qd_balance}).
$\hfill\blacksquare$

\vspace{2mm} \noindent \textbf{\textit{Proof of Theorem \ref{th_lb_bp}}}.
For any $k~(=k_1+k_2)$-tuple of columns in $\D$,
there are $s^{k_1}2^{k_2}$ level combinations in total.
Let $n=ts^{k_1}2^{k_2}+r_{n,k_1,k_2,s,2},
0\leq r_{n,k_1,k_2,s,2}<s^{k_1}2^{k_2}$.
According to (\ref{B_l1_lk}), to minimize $B_{l_1,\dots,l_k}$,
the frequencies of all the possible level combinations
in the column group
$\{\d_{l_1},\dots,\d_{l_k}\}$ should be as equal as possible.
Therefore, if there are $r_{n,k_1,k_2,s,2}$ level combinations
that occur $t+1$ times and $s^{k_1}2^{k_2}-r_{n,k_1,k_2,s,2}$
level combinations that appear $t$ times,
$B_{l_1,\dots,l_k}$ would reach its lower bound,
$r_{n,k_1,k_2,s,2}\left(1-\frac{r_{n,k_1,k_2,s,2}}{s^{k_1}2^{k_2}}\right)$.
By (\ref{Bk}), the lower bound of $B_k(\D)$ is obtained
immediately,
\begin{align}\label{lb_bk}
B_k(\D)\geq \sum_{\Omega}r_{n,k_1,k_2,s,2}
\left(1-\frac{r_{n,k_1,k_2,s,2}}{s^{k_1}2^{k_2}}\right)
\Bigg/ \dbinom{p+q}{k}.
\end{align}
By combining (\ref{qd_balance}) and (\ref{lb_bk}),
we have
\begin{align}\nonumber
 \text{QQD}^2(\D) \geq & -\left(\frac{5s+1}{4s}\right)^p
\left(\frac{4}{3}\right)^q+\left(\frac{5s+1}{4s}\right)^p
\left(\frac{11}{8}\right)^q
+\frac{1}{n^2}\left(\frac{5}{4}\right)^{p+q}
\sum_{k=1}^{p+q}\left(\frac{1}{5}\right)^k\\ \label{qd_lbs}
&\times\sum_{\Omega} r_{n,k_1,k_2,s,2}
\left(1-\frac{r_{n,k_1,k_2,s,2}}{s^{k_1}2^{k_2}}\right).
\end{align}
It's noted that \begin{align}\label{sum_omega}
\sum_{\Omega}=\sum_{k_1+k_2=k}\dbinom{p}{k_1}\dbinom{q}{k_2}.
\end{align}
The proof can be finished by substituting (\ref{sum_omega})
into (\ref{qd_lbs}).
$\hfill\blacksquare$

\end{document}